\documentclass[12pt]{article}

\usepackage{graphicx}
\usepackage{amscd}
\usepackage{theorem}

\newfont{\eulerfraktur}{eufm10 scaled\magstep1}

\newfont{\Bbb}{msbm10 scaled\magstep1}

\newcommand{\CC}{{\mbox{\Bbb C}}}

\newcommand{\PP}{{\mbox{\Bbb P}}}

\newcommand{\QQ}{{\mbox{\Bbb Q}}}

\def\bdm{\begin{displaymath}}
\def\be{\begin{equation}}

\def\CMF{\overline{{\cal M}}_g}

\def\edm{\end{displaymath}}
\def\ee{\end{equation}}

\def\lra{{\longrightarrow}}

\def\cocoa{{\hbox{\rm C\kern-.13em o\kern-.07em C\kern-.13em o\kern-.15em A}}}
\def\wh{\hat}
\def\wgt{{\rm wt}}
\renewcommand{\:}{\colon\,}
\newcommand{\x}{\times}
\newcommand{\ol}{\overline}

\newcommand{\w}{\omega}

\def\proof{{\bf Proof.}\,\,\,}

\def\qed{\hfill\vrule height4pt width4pt depth0pt}

\theoremstyle{change} 
\newtheorem{claim}{\hskip-0.2cm}[section]

\begin{document}

\title{Special ramification loci 
on the double product of a general curve
}

\author{{\small{{ C. CUMINO $^{(\natural)}$}, 
E. ESTEVES $^{(\flat)}$  AND  
{L. GATTO} $^{(\natural)}$ }}
\thanks{Work partially sponsored by MURST 
(Progetto Nazionale ``Geometria delle Variet\`a Proiettive'' 
Coordinatore Sandro Verra), and supported by 
GNSAGA-INDAM. The second author was also supported
by CNPq, Proc. 478625/03-0 and 301117/04-7, 
and CNPq/FAPERJ, Proc. E-26/171.174/2003.}\\
{ }\\
{\normalsize $^{(\natural)}$ Dipartimento di 
Matematica, Politecnico di Torino,}\\
{\normalsize Corso Duca degli Abruzzi 24, 
10129 Torino -- (ITALY)}\\{æ}
\\{\normalsize $^{(\flat)}$ Instituto Nacional de 
Matem\'atica Pura e Aplicada,}\\
{\normalsize Estrada Dona Castorina 110} \\ 
{\normalsize 22460-320 Rio de Janeiro RJ -- (BRAZIL)}}

\date{ }

\maketitle

\begin{abstract}
\noindent
Let $C$ be a general connected, smooth, projective 
curve of 
positive genus $g$. For each integer $i\geq 0$ we 
give formulas for the number of pairs 
$(P,Q)\in C\times C$ off the diagonal such that 
$(g+i-1)Q-(i+1)P$ is linearly equivalent to an 
effective divisor, and the number of pairs 
$(P,Q)\in C\times C$ off the diagonal such that 
$(g+i+1)Q-(i+1)P$ is linearly equivalent to a moving 
effective divisor.
\end{abstract}

\section{Introduction}\label{sec1}

Let $C$ be a general connected, smooth, projective 
curve of genus $g>0$. Put $C^2:=C\times C$, and let 
$\Delta\subset C^2$ be the diagonal. For each 
integer $i\geq 0$ consider the 
following loci on $C^2$:
     \begin{eqnarray*}
       D_i&:=&\{(P,Q)\in C^2-\Delta\,|\,
       h^0({\cal O}_C((g+i-1)Q-(i+1)P))>0\},\\
       E_i&:=&\{(P,Q)\in C^2-\Delta\,|\,
       h^0({\cal O}_C((g+i+1)Q-(i+1)P))>1\}.       
     \end{eqnarray*}
Our Proposition \ref{A2CC} 
claims that $D_i$ and $E_i$ are 
finite, and our main result, 
Theorem \ref{thm6DiazCuk}, gives 
formulas for the number of points in $D_i$ and $E_i$.

A formula for the number of points in $D_i$ appeared 
already as Lemma~6.3 on page 24 
of the seminal work by 
Diaz \cite{Diaz/Exc}, where the unnecessary extra 
hypotheses that $g\geq 2$ and $i\geq 2$ are made. 
Diaz used this formula to compute the class in 
the moduli space of genus-$g$ stable curves 
$\CMF$ of the closure $\ol{\cal D}_g$ 
of the locus of smooth 
curves $C$ having a Weierstrass point $P$ 
of type $g-1$, i.e. such that 
$h^0({\cal O}_C((g-1)P))\geq 2$.

Later on, Cukierman \cite{Cukie} gave a formula for 
the class in $\CMF$ of the closure $\ol{\cal E}_g$ 
of the locus of smooth 
curves $C$ containing a Weierstrass point $P$ of 
type $g+1$, 
i.e. such that $h^0({\cal O}_C((g+1)P))\geq 3$. 
He did not follow in Diaz's 
footsteps for this formula, but 
rather observed that the union 
$\ol{\cal D}_g\cup\ol{\cal E}_g$ is the branch 
locus of the Weierstrass divisor on the 
``universal'' curve over $\CMF$, and used a Hurwitz 
formula with singularities to compute the class of 
this branch divisor.
 
Had Cukierman followed in Diaz's footsteps, he 
would probably have found he needed a formula for 
the number of points in $E_i$. We give this formula 
here.

In fact, in a sense to be explained below, it is 
slightly easier to obtain the number of points in 
$E_i$ than in $D_i$, though we obtain both in a 
quite integrated form here. 
To obtain these numbers, the 
natural procedure is to use Porteous formula 
to compute the virtual classes of certain 
natural ramification schemes $D_i^+$ and $E_i^+$ of 
maps of vector bundles on $C^2$; see 
Subsection \ref{srs}. Set-theoretically, $D_i^+$ and 
$E_i^+$ are given exactly as $D_i$ and $E_i$, but 
without the restriction that the pair $(P,Q)$ lies 
off $\Delta$. 

The problem is that $D_i^+$ and $E_i^+$ are both 
larger than $D_i$ and $E_i$. Indeed, $E_i^+$ is the 
union of $E_i$ with the set of points $(P,P)$ such 
that $P$ is a Weierstrass point of $C$ and, worse, 
$D_i^+$ is the union of $D_i$ and the whole diagonal 
$\Delta$. Since $E_i^+$ is finite, 
Porteous formula does give an 
expression for the number of points in $E_i^+$, with 
weights, and thus at least an upper 
bound for the number of points in $E_i$. But it 
does not a priori give any information on $D_i$.

To compute the number of points in $D_i$ and $E_i$, 
we use the fact that, by the Riemann-Roch Theorem, 
the union of $D_i$ and $E_i$ is the locus $SW_i$ of 
pairs $(P,Q)\in C^2-\Delta$ such that $Q$ is a 
special ramification point of the complete linear 
system $H^0(\omega_C((i+1)P))$, where $\omega_C$ 
is the canonical bundle of $C$. 

We give $SW_i$ a scheme structure as follows. First, 
we consider the ramification divisor 
$Z_i\subset C^2$ of the 
family of linear systems $H^0(\omega_C((i+1)P))$ 
parameterized by $P\in C$. Our Proposition \ref{prop5} 
implies that 
$Z_i$ contains $\Delta$ with multiplicity exactly 
$g$. Set $W_i:=Z_i-g\Delta$. Our 
Proposition \ref{prop7} 
gives an expression for the cycle $[W_i]$, and 
our Proposition \ref{lem} claims that $W_i$ is 
nonsingular. Furthermore, in Subsection \ref{srs} we 
observe that the branch divisor of $W_i$ with 
respect to the projection $p_1\: C^2\to C$ over the 
first factor has support $SW_i$. We give 
$SW_i$ the structure of this branch divisor. 

The advantage of considering $SW_i$ is 
that it is quite easy to compute its degree. 
Indeed, $[SW_i]$ is the second Chern class 
of the bundle of first-order relative jets of $p_1$ 
with coefficients in ${\cal O}_{C^2}(W_i)$. 
Having an expression for $[W_i]$ we derive very 
quickly an expression for $\int_{C^2}[SW_i]$ in 
Proposition \ref{prop07}

Now, giving $D_i$ and $E_i$ the subscheme 
structures induced from $D_i^+$ and $E_i^+$, our 
Proposition \ref{A2CC} shows that, as 0-cycles, 
        $$[D_i]+[E_i]=[SW_i].$$
Actually, $D_i$ and $E_i$ are reduced. Indeed, in 
the proof of Theorem \ref{thm6DiazCuk} we show 
that the weight of $(P,Q)$ in $[SW_i]$ is at most 2, 
and the maximum weight is achieved if and only if 
$(P,Q)\in D_i\cap E_i$. 

Now, as we already know $\int_{C^2}[SW_i]$, 
it is enough to compute either $\int_{C^2}[D_i]$ or 
$\int_{C^2}[E_i]$. As mentioned above, we compute 
the latter. In fact, we can get $\int_{C^2}[E_i^+]$ 
using Porteous formula, and a local analysis, done 
in Proposition \ref{A2CC}, 
shows that the weight of $(P,P)$ 
in $[E_i^+]$ is equal to $g+1$ for each Weierstrass 
point of $C$. Thus $\int_{C^2}[E_i]$ follows.

In a second article \cite{CEG}, we 
show how the knowledge of the number of points in 
$E_i$ can be used to compute the class of 
$\ol{\cal E}_g$ in $\ol{\cal M}_g$. 
This computation is not 
straightforward as, following in Diaz's footsteps, 
we have to determine 
the limits of 
special Weierstrass points of type $g+1$ on stable 
curves with just one node. This is the main result 
of \cite{CEG}. 

The limits of special 
Weierstrass points of type $g-1$ were computed by 
Diaz, using admissible covers. However, the same 
method does not apply to points of type $g+1$. For 
those we apply in \cite{CEG} 
the theory of limit linear series in 
a rather new way, using 2-parameter families. 
Actually, as in the present article, we use an 
integrated approach in \cite{CEG} that 
yields simultaneously 
the limits of special Weierstrass points 
of both types, and also formulas for 
the classes of
both $\ol{\cal D}_g$ and $\ol{\cal E}_g$.

Here is a layout of the article. In 
Section \ref{sec2} we 
review the theory of linear systems and ramification 
on 
a
smooth curve $C$, introduce the linear 
systems 
we will consider in the remaining of the article, 
$H^0(\omega_C((i+1)P))$ for $P\in C$, and prove a 
preliminary result about them. In 
Section \ref{sec3}, 
assuming $C$ is general, we obtain through 
degeneration methods results that bound the order 
sequence of $H^0(\omega_C((i+1)P))$ at any point of 
$C$. In Section \ref{sec4}, 
we describe the structure 
of the ramification divisor $Z_i\subset C^2$ of the 
family of linear systems $H^0(\omega_C((i+1)P))$ 
parameterized by $P\in C$. 
Finally, in Section \ref{sec5} we define the loci 
$D_i$ and $E_i$ and compute their number of points, 
through the study of the locus $SW_i$ of special 
ramification points of the family 
$H^0(\omega_C((i+1)P))$ for $P\in C$.

We thank Nivaldo Medeiros for discussions on related 
topics. Also, we acknowledge the use of \cocoa 
\cite{cocoa} for some of the computations.

\section{Setup}\label{sec2}

\begin{claim} 
{\rm ({\it Linear systems and ramification})
Let $C$ be a {\it smooth curve}, that is, 
a projective, connected,
smooth scheme of dimension 1 over $\CC$. Denote by $\w_C$
its canonical sheaf. Let $g:=h^0(C,\w_C)$, 
the genus of $C$.

Let $V$ be a $\CC$-vector space of sections
of a line bundle $\cal L$ on $C$. We call $V$ a
{\it linear system}. The linear system is
called {\it complete} if $V=H^0(C,{\cal L})$.
Let $r:=\dim V-1$ and $d:=\deg{\cal L}$. We call
$r$ the {\it rank} of $V$ and $d$ its {\it degree}. We say
as well that $\dim V$ is the {\it dimension} of $V$.

For each point $P$ of $C$, and each integer $j\geq 0$,
let $V(-jP)$ denote the vector subspace of $V$ of sections
of $\cal L$ that vanish with order at least
$j$ at $P$. We say that $j$ is an {\it order} of
$V$ at $P$ if $V(-jP)\neq V(-(j+1)P)$. There are
$r+1$ orders, which, in an increasing sequence, will be
denoted by
       $$\epsilon_0(V,P),\,\epsilon_1(V,P),\,\dots,\,
       \epsilon_r(V,P).$$

For each integer $\ell\geq 0$ and 
each line bundle ${\cal M}$ on
$C$, let ${\cal J}^\ell_C({\cal M})$ 
be the bundle of {\it jets},
or {\it principal parts}, of $\cal M$ 
truncated in order $\ell$.
Consider the map of rank-$r$ bundles,
      $$V\otimes{\cal O}_C\longrightarrow
      {\cal J}^r_C({\cal L}),$$
locally obtained by differentiating up to 
order $r$ the sections
of $\cal L$ in $V$. The {\it wronskian} $w_V$ of 
$V$ is the (nonzero) section of
      $${\cal L}^{\otimes r+1}\otimes
      \w_C^{\otimes r(r+1)/2}$$
induced by taking determinants in the
above map of bundles.

For each point $P$ of $C$, the {\it weight} $\wgt_V(P)$ of
$P$ in $V$ is the order of vanishing of $w_V$ at $P$.
We call $P$ a {\it ramification point} of $V$ if
$\wgt_V(P)>0$; otherwise we call $P$ {\it ordinary}.
A local analysis yields the formula
       $$\wgt_V(P)=\sum_{j=0}^r
       (\epsilon_j(V,P)-j).$$
We call $P$ a {\it simple} ramification point if
$\wgt_V(P)=1$; otherwise we call $P$ {\it special}.
The point $P$ is special if and only if the section
       $$Dw_V\in H^0\Big(C,{\cal J}^1_C
       ({\cal L}^{\otimes r+1}\otimes
       \w_C^{\otimes r(r+1)/2})\Big),$$
locally obtained from $w_V$ by differentiating, 
vanishes at $P$.

The {\it total weight} of the ramification points of
$V$ is the (finite) sum
      $$\wgt_V:=\sum_{P\in C}\wgt_V(P).$$
It is equal to the degree of the line bundle of which
$w_V$ is a section, that is,
      $$\wgt_V=(r+1)(d+(g-1)r),$$
a formula usually referred to as the
{\it Brill--Segre} or {\it Pl\"ucker formula}.

The {\it canonical system} is the complete linear
system of sections of $\w_C$. Its rank is $g-1$, and its
degree is $2(g-1)$. For each point $P$ of $C$, its
{\it Weierstrass weight} $\wgt(P)$
is its weight in the canonical system, and the
{\it Weierstrass order sequence} at
$P$ is the increasing sequence of orders at $P$ of
the canonical system.

For each integer $i\geq -1$ and each $P\in C$,
let $V_C(i,P)$ denote the complete linear system of
sections of $\w_C((i+1)P)$.}
\end{claim}

\begin{claim}\label{prop5}{\bf Proposition.} 
Let $C$ be a smooth
curve of genus $g$. For each integer
$i\geq 0$ and each $P\in C$, the following two
statements hold for $V:=V_C(i,P)$:
\begin{enumerate}
\item\label{it5}
The weight $\wgt_V(P)$ of $P$ as a ramification point of
$V$ satisfies
      \be
      \wgt_V(P)=g+\wgt(P),\label{eq:15}
      \ee
where $\wgt(P)$ is the Weierstrass weight of $P$.
\item \label{it4} The total weight $\wgt_V$ of
the ramification points of $V$ satisfies
      \be
      \wgt_V=g(g+i)^2.\label{eq:14}
      \ee
\end{enumerate}
\end{claim}

\proof From the Riemann--Roch theorem, for
each $j=0,\dots,i$,
      $$\dim V(-jP)=g+i-j.$$
In particular, comparing dimensions, we get
      $$V(-iP)=V(-(i+1)P)=H^0(C,\w_C).$$
Hence, the order sequence of $V$ at $P$ is
      $$0,\, 1,\, \dots,\, i-1,\,
      (i+1)+\epsilon_0,\, (i+1)+\epsilon_1,\,
      \dots,\, (i+1)+\epsilon_{g-1},$$
where $\epsilon_0,\,\epsilon_1,\,\dots,\,\epsilon_{g-1}$ 
is the Weierstrass order sequence at $P$. Thus
      $$\wgt_V(P)=\sum_{k=0}^{g-1}(i+1+\epsilon_k-i-k)
      =g+\wgt(P).$$

The second statement is a direct application of
the Brill--Segre formula, using that the rank of $V$
is $g+i-1$ and its degree is $2g-2+(i+1)$.
\qed

\section{The general curve}\label{sec3}

\begin{claim}\label{WiPgen}{\bf Proposition.}
Fix an integer $i_0\geq 0$.
Let $C$ be a general smooth curve of genus
$g\geq 1$. Then the following two statements hold for
each nonnegative integer $i\leq i_0$:
\begin{enumerate}
\item For a general point $P$ of $C$, the linear system
$V_C(i,P)$ ramifies at $P$ with weight $g$, and has
otherwise at most simple ramification points.
\item For any two points $P$ and $R$ of $C$,
    \begin{equation}\label{WP-R1}
      h^0(C,\w_C((i+1)P-(g+i)R))\leq 1
    \end{equation}
\vskip-0.8cm
    \begin{equation}\label{WP-R2}
      h^0(C,\w_C((i+1)P-(g+i+2)R)=0.
    \end{equation}
\end{enumerate}
\end{claim}

\proof Let us first observe that the property required
of $C$ is open. Indeed, let $f\: X\to S$ be any family of
smooth curves, that is, a projective, smooth map with
connected fibers of dimension 1. Consider the fibered
product $X^{(2)}:=X\x_S X$ of two copies of $f$, and
denote by $p_1$ and $p_2$ the projection maps. Denote by
$\Delta$ the diagonal subscheme of $X^{(2)}$. Let
$\w_f$ denote the relative canonical bundle of $f$.
Then $\w:=p_2^*\w_f$ is the relative canonical bundle
of $p_1$. Let
        $${\cal V}:=p_{1*}(\w((i+1)\Delta)).$$
A fiberwise analysis shows that 
$\cal V$ is a bundle of rank $g+i$
with formation commuting with base change. 
For each integer
$\ell\geq 0$, denote by $\cal J^\ell$ the 
bundle 
of rank $\ell+1$
of $p_1$-relative
jets of $\w((i+1)\Delta)$ truncated in 
order $\ell$, and denote by
$\psi_\ell\: p_1^*{\cal V}\to{\cal J}^\ell$ 
the map locally
obtained by differentiating the sections 
of $\w((i+1)\Delta)$ up
to order $\ell$ along the fibers of $p_1$. 
Let $W_{i,1}$ (resp.
$W_{i,2}$) be the closed subset of $X^{(2)}$ 
where $\psi_{g+i-1}$
(resp. $\psi_{g+i+1}$) has rank at most 
$g+i-2$ (resp. $g+i-1$).
Also, let $W_i$ be the closed subset 
of $X^{(2)}$ where
$\psi_{g+i-1}$ has rank at most 
$g+i-1$.
By
Proposition~\ref{prop5}, $W_i$ 
contains $\Delta$ with multiplicity
$g$. Let $W'_i:=W_i-g\Delta$ and 
$Z_i:=\Delta\cap W'_i$. Let
$W''_i\subset W'_i$ be the ramification scheme of the map
$p_1|_{W'_i}$. Let $U_i\subseteq S$ be the intersection of
$S-f(p_1(W_{i,1}\cup W_{i,2}))$ with $f(X-p_1(W''_i\cup Z_i))$.
Since $p_1$ is proper, and $f$ is both proper and open, $U_i$ is
an open subscheme of $S$. Let $U:=U_0\cap\cdots\cap U_{i_0}$. The
formation of $U$ commutes with base change. Thus a fiberwise
analysis reveals that $U$ consists of the set of points $s\in S$
such that the proposition holds for $C:=X(s)$.

Now, keeping in mind the existence of a versal family of
smooth curves, it is enough to exhibit a single curve $C$
for which the statement holds.
We will actually show a somewhat 
stronger existence result:

\begin{claim}\label{WCQ}{\bf Lemma.}
Fix nonnegative integers $i_0$ and $j_0$.
Let $g$ be a positive integer.
Then there is a smooth
pointed curve $(C,Q)$ of genus $g$ for which
the following three statements hold for
each nonnegative integers $i\leq i_0$ and $j\leq j_0$:
\begin{enumerate}
\item The linear system $V_C(j,Q)$ ramifies at
$Q$ with weight g, and has otherwise at most simple
ramification points.
\item For each $P\in C$ distinct from $Q$, either $Q$
is an ordinary point or a simple ramification point of
$V_C(i,P)$.
\item For each $P\in C$ distinct from $Q$, the linear
system $V$ of sections of $\w_C((i+1)P+(j+1)Q)$ given by
        $$V:=H^0(\w_C((i+1)P))+H^0(\w_C((j+1)Q))$$
satisfies
        $$\dim V(-(g+i+j)R)\leq 1\quad{\rm and}\quad
        V(-(g+i+j+2)R)=0$$
for each $R\in C$ distinct from $P$ and $Q$.
\end{enumerate}
\end{claim}

We will first see how the lemma implies the
proposition. Set $j_0=i_0$, and consider the pointed curve
$(C,Q)$ given by the lemma. Then the two statements of
Proposition~\ref{WiPgen} hold for $C$. Indeed,
the first statement holds for $P=Q$, whence
for $P$ in a neighborhood of $Q$, that is,
for a general $P$.

As for the second statement, first notice that 
(\ref{WP-R1}) and (\ref{WP-R2}) 
hold for $P=Q$ and every $R\in C$, a consequence of
the first statement of the lemma for $j:=i$. They hold
as well for $R=Q$ and every $P\in C$ distinct from $Q$, a
consequence of the second statement of the lemma. 
Furthermore, they hold for $R=P$ and any $P\in C$. 
Indeed, the first statement of the lemma for $j:=0$ 
implies that the canonical linear system
has at most simple ramification points. Thus
$h^0(\w_C((1-g)P))\leq 1$ and $h^0(\w_C(-(g+1)P))=0$.

Finally, fix a point $P\in C$ distinct from $Q$ and a 
point $R\in C$ distinct from $P$ and $Q$. For $j:=0$, 
the linear system $V$ defined in the lemma is the 
system of sections of $\w_C((i+1)P+Q)$ that are zero 
on $Q$. Since $R\neq Q$, the third statement of the
lemma yields (\ref{WP-R1}) and (\ref{WP-R2}).

It is thus enough to prove the lemma, 
what we do below.\qed

\vskip0.2cm

\proof (Lemma \ref{WCQ})
We will do induction on $g$.
The initial step is taken care of below.

Let $C$ be any elliptic curve and 
$Q\in C$ any point. Then the ramification points
of the complete linear system of sections of 
$\w_C((j+1)Q)$
are simple. (These are the $(j+1)^2$ points $R$
for which $Q-R$ is $(j+1)$-torsion, what includes $Q$.)
In fact, it follows from the Riemann--Roch
theorem that every complete linear system has only simple
ramification points. Thus Statements~1 and 2 of the lemma
hold. Now, given $P\in C$ distinct from $Q$, since the
vector subspace $V$ of $H^0(\w_C((i+1)P+(j+1)Q))$ defined
in  Statement~3 has codimension 1, the order sequence of
$V$ at a point $R$ is obtained either from
        $$0,\,1,\,\dots,\,g+i+j-1,\,g+i+j\quad{\rm or}\quad
        0,\,1,\,\dots,\,g+i+j-1,\,g+i+j+1$$
by removing an order. In any case, there is at most one order of
$V$ at $R$ above $g+i+j-1$, that is $\dim V(-(g+i+j)R)\leq 1$, and
all orders are at most $g+i+j+1$, that is $V(-(g+i+j+2)R)=0$.

Assume from now on that $g>1$, and that the claim holds for
smaller genera and any integers $i_0$ and $j_0$. We will
employ a degeneration technique in order to apply the
induction hypothesis.

Let $(Y,A)$ and $(Z,B)$ be nonrational smooth pointed curves
of genera $g_Y$ and $g_Z$, with $g_Y+g_Z=g$. From the
induction hypothesis, we may assume that the statements
of the lemma hold for $(C,Q)$ replaced by $(Y,A)$ and
all nonnegative integers $i\leq i_0$ and
$j\leq g_Z+i_0+j_0+1$, and for
$(C,Q)$ replaced by $(Z,B)$ and all
nonnegative integers $i\leq i_0$ and $j\leq g_Y+i_0+j_0+1$.

Let $C_0$ be the curve of compact type 
that is the union of
$Y$, of $Z$, and of a chain of rational curves
$E_1,\dots,E_{n-1}$ connecting $A$ to 
$B$, where $n\geq 2$.
Our convention is that 
$E_1$ contains $A$ and $E_{n-1}$ contains $B$.
Let $v$ be any integer such that $0<v<n$, and let
$Q_0$ be any point of $E_v$ that is not a node of $C_0$.

Let $S:={\rm Spec}(\CC[[t]])$, and denote its special point
by 0 and generic point by $\eta$. Since there are no
obstructions to deforming pointed nodal curves, there
are a projective, flat map $f\: X\to S$ and a section
$\lambda\: S\to X$ of $f$ such that
$(X(0),\lambda(0))=(C_0,Q_0)$ and
$(X(\eta),\lambda(\eta))$ is a smooth pointed curve over
the field of formal Laurent series $\CC[[t]][1/t]$.

Let $C$ be the base extension of $X(\eta)$ to the algebraic
closure of $\CC[[t]][1/t]$. Set $Q:=\lambda(\eta)$. It is
enough to see that the statements of the lemma hold for
$(C,Q)$. Indeed, the argument is quite standard, and is
summarized below. Though the pointed curve $(C,Q)$ is not
defined over $\CC$, it is defined over a finitely generated
extension $L$ of $\QQ$. If the statements of the lemma hold
for $(C,Q)$, they also hold for the base extension of
$(C,Q)$ over any algebraically closed field containing $L$.
But, since $\CC$ has many transcendentals
over $\QQ$, there is an algebraically closed field
containing $L$ which is isomorphic to $\CC$. 
So, if the statements of the lemma hold for $(C,Q)$, 
they hold as well for some pointed curve over $\CC$.

Now, any finite set of points of $C$ is 
defined over a finite field extension of $\CC[[t]][1/t]$.
Replacing $S$ by its normalization in this field extension,
we may assume that these are rational points of $X(\eta)$,
and thus that there are sections of $f$ intersecting
$X(\eta)$ at them. By making a further base extension, if
necessary, and a sequence of blowups at the singular
points of the special fiber, we may assume that the total
space $X$ is regular, and that these sections factor through
the smooth locus of $f$. The compensation for this
is a change of the special fiber. However, the special fiber
will have the same specification as the $C_0$ we described
above. Thus, no confusion will ensue if we keep calling
by $C_0$ this new fiber. Also, the section $\lambda$ can be
extended to a section of this new family.

Now, let $P$ and $R$ be points of $C$ with $P$ distinct from
$Q$ and $R$ distinct from $P$ and $Q$. As we mentioned
above, we may assume there are sections $\gamma\:S\to X$ and
$\rho\:S\to X$ through the smooth locus of $f$ such that
$\gamma(\eta)=P$ and $\rho(\eta)=R$. Set $P_0:=\gamma(0)$
and $R_0:=\rho(0)$. Let $\Gamma$ and $\Lambda$ be the images
of $\gamma$ and $\lambda$, respectively.

Fix nonnegative integers $i\leq i_0$ and $j\leq j_0$.
Let $\w$ be the relative dualizing bundle of $f\: X\to S$.
Let $V_\eta$ be the linear system of sections of the line
bundle $\w(\eta)((i+1)P+(j+1)Q)$ given by
       $$V_\eta:=H^0(\w(\eta)((i+1)P))+
       H^0(\w(\eta)((j+1)Q)).$$
Assume that $R$ is a ramification point of $V_\eta$. To
prove the statements of the lemma hold for $(C,Q)$, it is
enough to prove the following three statements:
\begin{enumerate}
    \item For $i=0$, the system $V_\eta$ ramifies at $Q$ with weight
    $g$, and $R$ is a simple ramification point of $V_\eta$.
    \item For $j=0$, the point $Q$ is a ramification point of
    $V_\eta$ of weight $g+i$ or $g+i+1$.
    \item $\dim V_\eta(-(g+i+j)R)\leq 1$ and
    $V_\eta(-(g+i+j+2)R)=0$.
\end{enumerate}

We will employ techniques of limit linear series, 
from
\cite{MathCont}, to show the
above three statements. 
There are two cases to consider:

\vskip0.2cm

{\noindent\emph{Case 1:} Assume that $P_0\in E_u$
for some $u$.}

\vskip0.1cm

Since $C_0$ is of compact type, there is an effective
divisor $D$ of $X$ supported on $C_0$ such that, letting
       $${\cal L}:=\w((i+1)\Gamma+(j+1)\Lambda+D),$$
we have ${\cal L}|_{E_m}\cong{\cal O}_{E_m}$ for each
$m=1,\dots,n-1$,
       $${\cal L}|_Z\cong\w_Z((g_Y+i+j+3)B)\quad
       {\rm and}\quad
       {\cal L}|_Y\cong\w_Y((1-g_Y)A).$$
Since, from the induction hypothesis, $A$ is a ramification
point of weight $g_Y$ of the complete linear system of
sections of $\w_Y(A)$, the point $A$ is not a Weierstrass
point of $Y$. Then $V:=H^0(X,{\cal L})\cap V_\eta$ restricts
to a linear system $V_Z$ of dimension $g+i+j$ of sections of
$\w_Z((g_Y+i+j+3)B)$. Also from the induction hypothesis,
$B$ is not a Weierstrass point of $Z$. So the order
sequence of $B$ in the complete linear system of sections of
$\w_Z((g_Y+i+j+3)B)$ is
       $$0,\,1,\,\dots,\,g_Y+i+j+1,\,g_Y+i+j+3,\,\dots,\,
       g+i+j+2.$$
As a consequence, the weight $w_B$ of $B$ as a ramification
point of the linear system $V_Z$ satisfies
       \be\label{wB1}
       w_B\leq 2(g_Y+i+j)+3g_Z,
       \ee
with equality if and only if $V_Z=H^0(\w_Z((g_Y+i+j+1)B))$.

Analogously, choosing an appropriate $D$, we obtain a linear
system $V_Y$ of dimension $g+i+j$ of sections of
$\w_Y((g_Z+i+j+3)A)$, and the weight $w_A$ of $A$
as a ramification point of $V_Y$ satisfies
       \be\label{wA1}
       w_A\leq 2(g_Z+i+j)+3g_Y,
       \ee
with equality if and only if $V_Y=H^0(\w_Y((g_Z+i+j+1)A))$.

Let $r:=g+i+j-1$. Using the Pl\"ucker formula, the number
$N$ of ramification points of $V_Y$ and $V_Z$ on
$(Y-A)\cup(Z-B)$, counted with their respective weights,
satisfies
       \begin{eqnarray*}
       N&=&(r+1)\Big((2g_Z+g_Y+i+j+1)+r(g_Z-1)\Big)-w_B\\
       &+&(r+1)\Big((2g_Y+g_Z+i+j+1)+r(g_Y-1)\Big)-w_A\\
       &=&N'+5g+4(i+j)-w_A-w_B,
       \end{eqnarray*}
where
       $$N':=(r+1)\Big((2g+i+j)+r(g-1)\Big)-2g-i-j.$$

Now, from the theory of limit linear series, each one of the
ramification points of $V_Y$ or $V_Z$ on $(Y-A)\cup(Z-B)$
is a limit of ramification points of $V_\eta$, and its
weight as a ramification point is the sum of the weights of
the ramification points of $V_\eta$ converging to it.
Besides those, since $P$ and $Q$ are ramification points
of $V_\eta$ with weights at least $g+j$ and $g+i$,
respectively, the points $P_0$ and $Q_0$ appear
as limits of ramification points of $V_\eta$ with weights
summing up to at least $2g+i+j$. Thus, from the Pl\"ucker
formula, at most $N'$ ramification points of $V_\eta$,
counted with their weights, converge to $(Y-A)\cup(Z-B)$. So
       $$5g+4(i+j)-w_A-w_B\leq 0.$$

However, Inequalities (\ref{wB1}) and (\ref{wA1}) for $w_B$
and $w_A$ yield the opposite inequality:
       $$5g+4(i+j)-w_A-w_B\geq 0.$$
Thus, equalities hold, and hence
       $$V_Y=H^0(\w_Y((g_Z+i+j+1)A))\quad{\rm and}\quad
       V_Z=H^0(\w_Z((g_Y+i+j+1)B)).$$
In addition, $P$ and $Q$ are ramification points of
$V_\eta$ of weights $g+j$ and $g+i$, respectively, and
all the other ramification points of $V_\eta$ converge
to $(Y-A)\cup(Z-B)$. In particular, Statement 2 
and the first part of Statement~1 are shown.

Now, since $R$ is a ramification point of $V_\eta$, and
$R$ is distinct from $P$ and $Q$, we have
$R_0\in(Y-A)\cup(Z-B)$. So $R_0$ is a
ramification point of either $V_Y$ or $V_Z$. From the
induction hypothesis, the complete linear systems of
sections of $\w_Y((g_Z+i+j+1)A)$ and $\w_Z((g_Y+i+j+1)B)$
have at most simple ramification points, other than $A$ or
$B$. Thus $R$ is the unique ramification point of $V_\eta$
converging to $R_0$ and its weight is 1.
So the remainder of Statement 1 is shown.

As for Statement 3, assume, 
without loss of generality, that
$R_0\in Z$. Set
$n:=\dim V_\eta(-(g+i+j)R)$, and let
$\sigma_1,\dots,\sigma_n$ form a $\CC[[t]]$-basis of
$V\cap V_\eta(-(g+i+j)R)$. Their restrictions to $Z$ are
sections of $V_Z$ vanishing with multiplicity at least
$g+i+j$ on $R_0$. Assume, by contradiction, that
$n\geq 2$. Since $R_0$ is a simple ramification point of
$V_Z$, the sections $\sigma_1|_Z,\dots,\sigma_n|_Z$
are linearly dependent. Thus, there is
a nonzero $n$-tuple $(c_1,\dots,c_n)\in\CC^n$ such that
$c_1\sigma_1+\cdots+c_n\sigma_n$ vanishes on $Z$, and hence
on the whole $C_0$. Thus
       \be\label{sigma}
       c_1\sigma_1+\cdots+c_n\sigma_n=t\sigma
       \ee
for some $\sigma\in H^0(X,{\cal L})$. Also $\sigma\in
V_\eta(-(g+i+j)R)$, and hence $\sigma$ is a $\CC[[t]]$-linear
combination of $\sigma_1,\dots,\sigma_n$. Plugging this linear
combination in (\ref{sigma}) we obtain a nontrivial
$\CC[[t]]$-linear relation among the sections $\sigma_i$, a
contradiction. Thus $n\leq 1$. A similar analysis, using that
$V_Z(-(g+i+j+2)R_0)=0$, shows that $V_\eta(-(g+i+2)R)=0$,
finishing the proof of Statement 3.

\vskip0.2cm

{\noindent\emph{Case 2:} Assume $P_0$ belongs to either $Y$
or $Z$.}

\vskip0.1cm

Without loss of generality, we may assume that $P_0\in Z$.
Again, since $C_0$ is of compact type, there is an effective
divisor $D$ of $X$ supported on $C_0$ such that, letting
       $${\cal L}:=\w((i+1)\Gamma+(j+1)\Lambda+D),$$
we have ${\cal L}|_{E_m}\cong{\cal O}_{E_m}$ for each
$m=1,\dots,n-1$,
       $${\cal L}|_Z\cong\w_Z((i+1)P_0+(g_Y+j+2)B)
       \quad{\rm and}\quad
       {\cal L}|_Y\cong\w_Y((1-g_Y)A).$$
As before, $V:=H^0(X,{\cal L})\cap V_\eta$ restricts to a
linear system $V_Z$ of dimension $g+i+j$ of sections of
$\w_Z((i+1)P_0+(g_Y+j+2)B)$.

Now,
       $$V\supseteq H^0(X,\w((i+1)\Gamma))+
       H^0(X,\w((j+1)\Lambda+D)).$$
Reasoning as in Case 1, we can show that
$H^0(X,\w((j+1)\Lambda+D))$ restricts to
$H^0(\w_Z((g_Y+j+1)B))$.
On the other hand, the exact sequence
       $$0\to H^0(\w_Z((i+1)P_0))\to
       H^0(\w((i+1)\Gamma)|_{C_0})
       \to H^0(\w_Y(A))$$
shows that $h^0(\w((i+1)\Gamma)|_{C_0})=g+i$, and hence
that $H^0(X,\w((i+1)\Gamma))$ restricts to a vector
subspace of $H^0(\w_Z((i+1)P_0+B))$ containing
the subspace $H^0(\w_Z((i+1)P_0))$. Thus
       $$V_Z\supseteq H^0(\w_Z((i+1)P_0))+
       H^0(\w_Z((g_Y+j+1)B)),$$
and a dimension count shows that equality holds.

The weight $w_B$ of $B$ as a ramification point of $V_Z$
depends on its weight as a ramification point of
$V_Z(i,P_0)$. Now, from the induction hypothesis,
$B$ is either an ordinary point or a simple ramification
point of $V_Z(i,P_0)$.
Hence, the order sequence at $B$ of the linear system
$V_Z$ is either
       $$1,\,2,\,\dots,\,g_Y+j,\,g_Y+j+2,\,g_Y+j+3,\,
       \dots,\,g+i+j,\,g+i+j+1$$
or
       $$1,\,2,\,\dots,\,g_Y+j,\,g_Y+j+2,\,g_Y+j+3,\,
       \dots,\,g+i+j,\,g+i+j+2.$$
At any rate,
       \begin{equation}\label{wB2}
       w_B\leq g_Y+j+2(g_Z+i)+1.
       \end{equation}
Notice that, if $i=0$, then $B$ is an ordinary point of
$V_Z(0,P_0)$, as it is an ordinary point of $Z$,
and thus Inequality (\ref{wB2}) is strict.

On the other hand, let $D'$ be an effective divisor of $X$
supported on $C_0$ such that, letting
       $${\cal M}:=\w((i+1)\Gamma+(j+1)\Lambda+D'),$$
we have ${\cal M}|_{E_m}\cong{\cal O}_{E_m}$ for each
$m=1,\dots,n-1$,
       $${\cal M}|_Y\cong\w_Y((g_Z+i+j+3)A)
       \quad{\rm and}\quad
       {\cal M}|_Z\cong\w_Z((i+1)P_0-(g_Z+i)B).$$
Since, as mentioned above, $B$ is either an ordinary point
or a simple ramification point of $V_Z(i,P_0)$, we have that
$H^0(X,{\cal M})\cap V_\eta$ restricts to a linear
system $V_Y$ of dimension $g+i+j$ of sections of
$\w_Y((g_Z+i+j+3)A)$.

Since $A$ is not a Weierstrass point of $Y$, the sequence of
orders at $A$ of the complete linear system
of sections of $\w_Y((g_Z+i+j+3)A)$ is
       $$0,\,1,\,\dots,\,g_Z+i+j+1,\,g_Z+i+j+3,\,
       g_Z+i+j+4,\,\dots,\,g+i+j+2.$$
Since $V_Y$ has codimension 2 in $H^0(\w_Y((g_Z+i+j+3)A))$,
the weight $w_A$ of $V_Y$ at $A$ satisfies
       \begin{equation}\label{wA2}
       w_A\leq 2(g_Z+i+j)+3g_Y,
       \end{equation}
with equality if and only if $V_Y=H^0(\w_Y((g_Z+i+j+1)A))$.

As in Case 1, using the Pl\"ucker formula, the number $N$ of
ramification points of $V_Y$ and $V_Z$ on $(Y-A)\cup(Z-B)$,
counted with their respective weights, satisfies
       $$N=N'+4g+4i+3j-w_A-w_B,$$
where
       $$N':=(g+i+j)(2g+i+j)+(g+i+j)(g+i+j-1)(g-1)-g-i.$$
As in Case 1, since $Q$ is a ramification point of $V_\eta$
with weight at least $g+i$, there are at most $N'$
ramification points of $V_\eta$, counted with their
respective weights, converging to $(Y-A)\cup(Z-B)$. So
       $$4g+4i+3j-w_A-w_B\leq 0.$$
On the other hand, Inequalities (\ref{wB2}) and (\ref{wA2})
yield
       $$w_A+w_B\leq 4g+4i+3j+1.$$
In particular, $w_A\geq 2(g_Z+i+j)+3g_Y-1$, whence
       $$V_Y\subset H^0(\w_Y((g_Z+i+j+2)A)).$$
Also, $Q$ has weight $g+i$ or $g+i+1$ in $V_\eta$.
Thus Statement 2 is shown.
Furthermore, 
if $i=0$ we have $w_A+w_B=4g+4i+3j$. In this case,
$V_Y=H^0(\w_Y((g_Z+i+j+1)A))$ and 
$Q$ has weight $g+i$ in $V_\eta$, showing the first 
part of Statement 1.

If $Q$ has weight $g+i+1$ in $V_\eta$, all other ramification
points converge to $(Y-A)\cup(Z-B)$. If $Q$ has weight $g+i$,
there is at most one ramification point of 
$V_\eta$, other than $Q$,
converging outside $(Y-A)\cup(Z-B)$, and that point is 
simple. If $R$ is that point, then
$\dim V_\eta(-(g+i+j)R)\leq 1$ and 
$V_\eta(-(g+i+j+2)R)=0$ because
of the simplicity of $R$.

Assume now that $R_0\in (Y-A)\cup(Z-B)$. Let us first
consider the case $R_0\in Y-A$. In this case, since,
from the induction hypothesis,
the complete linear system of sections of
$\w_Y((g_Z+i+j+2)A)$ has at most simple
ramification points, other than $A$, we have
       \begin{eqnarray*}
       h^0(\w_Y((g_Z+i+j+2)A-(g+i+j)R_0)&\leq& 1,\\
       h^0(\w_Y((g_Z+i+j+2)A-(g+i+j+2)R_0)&=& 0.
       \end{eqnarray*}
Thus $\dim V_Y(-(g+i+j)R_0)\leq 1$ and $V_Y(-(g+i+j+2)R_0)=0$ as
well. It follows, as in Case 1, that 
$$
\dim V_\eta(-(g+i+j)R)\leq 1\quad {\rm and}\quad 
V_\eta(-(g+i+j+2)R)=0.
$$ 
Furthermore,
if $i=0$, since in this case $V_Y=H^0(\w_Y((g_Z+i+j+1)A))$, all
the ramification points of $V_Y$ distinct from $A$ are simple.
Thus $R_0$ is simple in $V_Y$, and hence $R$ is simple in
$V_\eta$.

Assume now that $R_0\in Z-B$. There are two cases to consider.
First, assume $R_0=P_0$. Since, by induction hypothesis, the
complete linear system of sections of $\w_Z((g_Y+j+1)B)$ has at
most simple ramification points other than $B$, the weight of
$P_0$ as a ramification point of $V_Z$ is either $g+j$ or $g+j+1$.
Since $P$ has at least weight $g+j$ in $V_\eta$, and $R\neq P$,
the latter must hold, and $R$ must be a simple ramification point
of $V_\eta$. In particular, $\dim V_\eta(-(g+i+j)R)\leq 1$ and
$V_\eta(-(g+i+j+2)R)=0$.

Finally, assume $R_0\neq P_0$. Then 
$$
\dim V_Z(-(g+i+j)R_0)\leq 1\quad {\rm and} \quad 
V_Z(-(g+i+j+2)R_0)=0
$$
from the induction hypothesis, and
hence $\dim V_\eta(-(g+i+j)R)\leq 1$ and $V_\eta(-(g+i+j+2)R)=0$.
Thus Statement~3 is shown. Also, if $i=0$, then
$V_Z=H^0(\w_Z((g_Y+j+1)B))$, and, since $R_0\neq P_0$, the weight
of $R_0$ in $V_Z$ is equal to its weight in the complete linear
system of sections of $\w_Z((g_Y+j+1)B)$. By induction hypothesis,
this weight is one, and thus $R$ is a simple ramification point of
$V_\eta$. So Statement~1 is shown.\qed

\begin{claim}\label{ws}{\bf Corollary.} 
If $C$ is a general
smooth curve of genus $g\geq 1$, then all its
Weierstrass points are simple.
\end{claim}

\proof 
Apply Statement 1 of
Proposition~\ref{WiPgen} for $i_0:=0$ and $i:=0$. 
\qed

\begin{claim}\label{arbD}{\bf Proposition.} 
Fix an integer $i_0\geq 0$. 
Let $C$ be a general smooth curve of genus $g\geq 1$. 
Then for any two distinct points $P$ and $R$ of $C$, 
and any nonnegative integer $i\leq i_0$,
       $$h^0(C,\omega_C((i+1)P-(g+i-2)R))=2.$$
\end{claim}

\proof A line bundle of degree 2 on an elliptic curve has 
(at most) 2 linearly independent sections. Thus we 
may assume $g\geq 2$. Also, for $i=0$, 
       $$h^0(C,\omega_C((i+1)P-(g+i-2)R)=
         h^0(C,\omega_C(-(g-2)R))=2,$$
since $R$ is at most a simple Weierstrass point of 
$C$, a consequence of 
Corollary~\ref{ws}. 
So we need 
only show the stated equality for integers $i>0$.

For each integer $j\geq 2$ (resp. $j\geq 1$), let 
$M_j$ be the moduli space of smooth curves (resp. 
let $M_{j,1}$ be the moduli space of smooth pointed 
curves) of genus $j$. 
Let $\ol M_j$ and $\ol M_{j,1}$ denote their 
respective compactifications by stable (resp. stable, pointed) curves. 
For each 
positive integer $i\leq i_0$, let 
$D^{(i)}\subseteq M_{g+i}$ be the subset 
parameterizing curves admitting a covering 
of degree at most $g+i-2$ of the projective line 
totally ramified at a point. 
By \cite{A}, Thm. 3.11, p. 333, the subvariety 
$D^{(i)}$ is irreducible of codimension 2. Let 
$\ol D^{(i)}\subset\ol M_{g+i}$ be the closure of 
$D^{(i)}$. 

Let $\mu_i\:M_{g,1}\times M_{i,1}\to 
\ol M_{g+i}$ be the natural map, associating to a 
pair of smooth pointed curves the stable uninodal curve 
which is the union of these curves identified at the 
marked points. Let $E^{(i)}:=\mu_i^{-1}(\ol D^{(i)})$. 
Let 
$\rho_i\:E^{(i)}\to M_g$ be the natural map, 
forgetting the second 
pointed curve and the marked point on the first curve. 
Since $C$ is general, we may assume that, 
for each $i=1,\dots,i_0$, the 
curve $C$ is parameterized 
by a point of $M_g$ over which the fiber of 
$\rho_i$ has minimum dimension. We claim this 
dimension is at most $3i-3$, whence less than 
$\dim M_{i,1}$. Indeed, if the dimension were 
larger, then $E^{(i)}$ would have codimension at 
most 1 
in $M_{g,1}\times M_{i,1}$, and hence would dominate 
$\ol D^{(i)}$ under $\mu_i$. So $\ol D^{(i)}$ would 
be contained in the boundary $\ol M_{g+i}-M_{g+i}$, 
an absurd. From the claim, for each $i=1,\dots,i_0$, 
the general smooth pointed curve $(Y_i,B_i)$ of genus 
$i$ is such that,
for any $P\in C$, the 
pair of pointed curves $((C,P),(Y_i,B_i))$ is 
not parameterized by $E^{(i)}$. Consequently, the 
stable uninodal curve $X_i$, union of $C$ and $Y_i$ 
with $P$ and $B_i$ identified, is 
parameterized by a point of 
$\ol M_{g+i}$
off $\ol D^{(i)}$, 
for each $i=1,\dots,i_0$. 

Suppose, by contradiction, that for certain 
distinct points $P$ and $Q$ of $C$, and a certain 
positive integer $i\leq i_0$, we have
       $$h^0(C,\w_C((i+1)P-(g+i-2)Q)\geq 3.$$
Put $g':=g+i$. Since, by Riemann--Roch, 
$h^0(C,\w_C((i+1)P-iQ))=g$, 
there is an integer $j$ with $2\leq j<g$ 
such that
       $$h^0(C,\w_C((i+1)P-(g'-j)Q)=
         h^0(C,\w_C((i+1)P-(g'-j-1)Q)=j+1.$$
Again by Riemann--Roch, 
       \begin{equation}
	 h^0(C,{\cal O}_C((g'-j)Q-(i+1)P))>
	 h^0(C,{\cal O}_C((g'-j-1)Q-(i+1)P)).
       \end{equation}
Thus, there is a map 
$\phi\:C\lra \PP^1$ of degree $g'-j$ such that 
$\phi^*(0)=(g'-j)Q$ and $\phi^*(\infty)\geq(i+1)P$. 
Let $i'$ be the integer such that $i'+1$ is the 
multiplicity of $P$ in 
$\phi^*(\infty)$. 
Then $i'\geq i$.

Set $Y:=Y_i$ and $B:=B_i$. 
Since $B$ is general, $B$ is not a Weierstrass 
point of $Y$. Thus, since $i'\geq i$, we have 
$h^0(Y,{\cal O}_Y(i'B))<h^0(Y,{\cal O}_Y((i'+1)B))$.  
So, there is a map 
$\psi\:Y\lra\PP^1$ of degree $i'+1$ such that 
$\psi^*(\infty)=(i'+1)B$. 

Putting together the maps 
$\phi$ and $\psi$, we may construct the 
covering with source $X_i$ depicted  
in Figure~1 below,
\begin{figure}[ht]
\begin{center}
\includegraphics[height=5cm, angle=0]{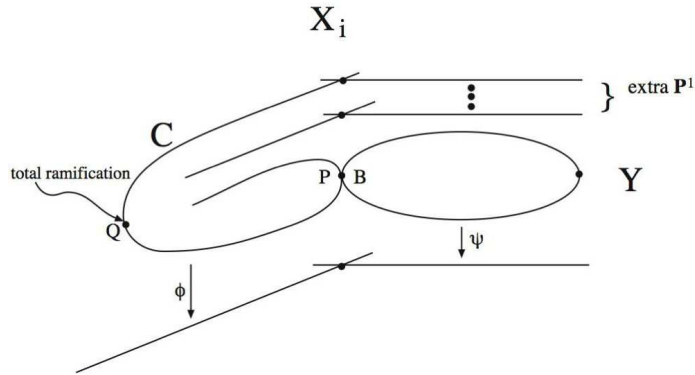}
\caption{The covering.}
\end{center}
\end{figure}
which can be represented by a point $[X_i]$ of the 
(compactification of the) Hurwitz scheme 
parameterizing (pseudo)admissible 
coverings of the projective line of 
degree $(g'-j)$ totally ramified at a 
point; see Remark \ref{rmk1}.
Since 
coverings of $\PP^1$ form a dense open subscheme of 
this compactification, the curve 
$X_i$ is limit of smooth curves equipped with a 
degree-$(g'-j)$ map to the projective line totally 
ramified at a 
point. Since $j\geq 2$, it follows that $[X_i]$ 
lies on the boundary of $D^{(i)}$,
a contradiction.\qed

\begin{claim}\label{rmk1}{\bf Remark.} {\rm The 
Hurwitz scheme we used in the proof of 
Proposition \ref{arbD} is 
mentioned in \cite{Diaz/Exc}, Section 5. It can be 
constructed following the same reasoning used in the 
construction of the Hurwitz scheme of (simple) 
admissible coverings, given in the proof 
of \cite{HaMu}, Thm. 4, p. 58. Also, the local 
descriptions of both schemes are the same, given on 
\cite{HaMu}, p. 62. From this description we see that 
the Hurwitz scheme is equidimensional. Now, there 
is a natural forgetful map from the Hurwitz scheme 
to a corresponding moduli space of pointed genus-0 
curves, taking a covering to its target. 
This map is finite and surjective, also by 
\cite{HaMu}, Thm. 4, p.~58. Since
the moduli spaces of pointed genus-$0$ curves 
are irreducible (see \cite{Kn} or \cite{Ke}),
it follows that each irreducible 
component of the Hurwitz scheme covers the target. So 
coverings of $\PP^1$ form a dense open subscheme of 
the Hurwitz scheme, a fact used in the proof 
of Proposition \ref{arbD}.}
\end{claim}

\begin{claim}\label{rmk2}{\bf Remark.} {\rm We tried 
to prove Proposition \ref{arbD} using the same induction 
argument used in the proof of Lemma \ref{WCQ}. 
However, we could not prove the 
initial step, that is, the following statement:}
{\it Let $C$ be a general elliptic 
curve, $Q\in C$ a general point, and $P\in C-\{Q\}$ any 
point. Let $i$ and $j$ be nonnegative integers. Then the 
linear system $V$ of sections of the line bundle 
$\w_C((i+1)P+(j+1)Q)$ 
generated by $H^0(\w_C((i+1)P))$ and $H^0(\w_C((j+1)Q))$ 
satisfies $\dim V(-(i+j-1)R)=2$ for 
each $R\in C-\{P,Q\}$.}
\end{claim}

\section{Weierstrass divisors}\label{sec4}

\begin{claim}\label{wronski} {\rm ({\it Wronski maps})
Let $C$ be a smooth curve of genus $g$. 
For each integer $j\geq 0$, consider the family of
linear systems $V_C(j,P)$ for $P$ varying on $C$.
More precisely, let
$p_1$ and $p_2$ denote the projections of $C\times C$ onto
the first and second factors, and 
$\Delta\subset C\times C$ the diagonal.
The relative canonical bundle of $p_1$
is simply the pullback $p_2^*\w_C$ of the canonical bundle
$\w_C$ of $C$. For each integer $j\geq 0$, let
       $${\cal L}_j:=p_2^*\w_C((j+1)\Delta),\quad
       {\cal E}_j:=p_{1*}{\cal L}_j.$$
Notice that, for each point $P$ of $C$,
identifying $\{P\}\times C$ with $C$
in the natural way,
${\cal L}_j|_{\{P\}\times C}=\w_C((j+1)P)$.
Also, as $h^0(\w_C((j+1)P))=g+j$ for
every $P\in C$, the sheaf ${\cal E}_j$ is a
bundle of rank
$g+j$ and ${\cal E}_j|_P=H^0(\w_C((j+1)P))$.

For each integer $\ell\geq 0$ and each line bundle
${\cal M}$ on $C\times C$, let
${\cal J}^\ell_{p_1}({\cal M})$
be the bundle 
of rank $\ell+1$ 
of $p_1$-relative jets of ${\cal M}$
truncated in order $\ell$. Let
       $$\rho_{j,\ell}\:p_1^*{\cal E}_j\to
       {\cal J}^\ell_{p_1}({\cal L}_j)$$
be the map of bundles locally obtained by 
differentiating 
up to order $\ell$ along the fibers of $p_1$
the sections of ${\cal L}_j$. We call $\rho_{j,\ell}$ 
a {\it Wronski map}.

The map $\rho_{j,g+j-1}$ is a map of bundles of the
same rank. Taking determinants, we get a section $z_j$
of the line bundle
        $$\bigwedge^{g+j}{\cal J}^{g+j-1}_{p_1}({\cal L}_j)
        \otimes\bigwedge^{g+j}p_1^*{\cal E}_j^{\vee},$$
which is naturally isomorphic, using the truncation
sequence of the bundles of jets, to
        $$p_2^*\w_C((j+1)\Delta)^{\otimes g+j}\otimes
        p_2^*\w_C^{\otimes (g+j)(g+j-1)/2}\otimes
        \bigwedge^{g+j}p_1^*{\cal E}_j^{\vee},$$
or more simply to
        $$p_2^*\w_C^{\otimes (g+j)(g+j+1)/2}
        \Big((g+j)(j+1)\Delta\Big)\otimes\bigwedge^{g+j}
        p_1^*{\cal E}_j^{\vee}.$$}
\end{claim}

\begin{claim}\label{jWd} 
{\rm ({\it Weierstrass divisors.}) 
Keep the notation used in Subsection \ref{wronski}.
Let $Z_j\subseteq C\times C$ denote the zero scheme 
of $z_j$. The section $z_j$ is a relative wronskian.
More precisely, for each $P\in C$, on $\{P\}\times C$,
identified with $C$ in the natural way,
the section $z_j$ restricts to the
wronskian of the linear system $V_C(j,P)$. Hence, $Z_j$
consists of the pairs $(P,Q)\in C\times C$ such that
$V_C(j,P)$ ramifies at $Q$. Now, since $z_j$ is nonzero,
being so on each fiber, $Z_j$ is a Cartier divisor.
By Proposition~\ref{prop5}, the divisor $Z_j$ intersects
each fiber $\{P\}\times C$ at $(P,P)$ with multiplicity
$g+\wgt(P)$, where $\wgt(P)$ is the weight of $P$ in the
canonical system of $C$. Thus $Z_j$ contains $\Delta$ with
multiplicity exactly $g$. Let
        $$W_j:=Z_j-g\Delta.$$
Then $W_j$ is, set-theoretically, the locus of pairs
$(P,Q)\in C\times C$ such that either $P=Q$ and $P$ is a
Weierstrass point of $C$, or $P\neq Q$ and $Q$ is a
ramification point of $V_C(j,P)$. We call $W_j$ the 
{\it $j$-th Weierstrass divisor of $C$.}}
\end{claim}

\begin{claim}\label{prop7}{\bf Proposition.} 
Let $C$ be a smooth curve of
genus $g\geq 1$ and $j$ a nonnegative integer. Let 
$\Delta$ be the diagonal of $C\times C$, and 
$p_1$ and $p_2$ the projections of $C\times C$ 
onto the indicated factors. Let 
$\w_C$ be the canonical bundle of $C$, and set 
$K_\ell:=c_1(p_\ell^*\w_C)$ for $\ell=1,2$. 
Let $W_j\subseteq C\times C$ be 
the $j$-th Weierstrass divisor of $C$. 
Then its class $[W_j]$ in the Chow group of
$C\times C$ satisfies
        \be
        [W_j]=\frac{1}{2}(g+j)(g+j+1)K_2+j(g+j+1)[\Delta]
        +\frac{1}{2}j(j+1)K_1.\label{eq:17}
        \ee
\end{claim}

\proof Use 
the notation in Subsections \ref{wronski} 
and \ref{jWd}. Since $W_j=Z_j-g\Delta$, and $Z_j$ is the
zero scheme of a section of the line bundle
        $$p_2^*\w_C^{\otimes (g+j)(g+j+1)/2}
        \Big((g+j)(j+1)\Delta\Big)\otimes\bigwedge^{g+j}
        p_1^*{\cal E}_j^{\vee},$$
we get
        \be
    [W_j]=\frac{1}{2}(g+j)(g+j+1)K_2+j(g+j+1)[\Delta]-
        p_1^*c_1({\cal E}_j).\label{eq:67}
        \ee
To finish, we need only show that
        \be
        c_1({\cal E}_j)=-\frac{1}{2}j(j+1)c_1(\w_C).
    \label{eq:c1Ei}
    \ee

We show (\ref{eq:c1Ei}) by induction on $j$.
First of all,
        $${\cal E}_0=p_{1*}p_2^*\w_C=
        H^0(\w_C)\otimes{\cal O}_C.$$
Since ${\cal E}_0$ is free, $c_1({\cal E}_0)=0$.

Assume now that $j>0$ and
$c_1({\cal E}_{j-1})=-(j(j-1)/2)c_1(\w_C)$.
Consider the natural short exact sequence
        $$0\to p_2^*\w_C(j\Delta)\to p_2^*\w_C((j+1)\Delta)
        \to p_2^*\w_C((j+1)\Delta)|_{\Delta}\to 0.$$
Since $H^1(\w_C(jP))=0$ for each $P\in C$,
applying $p_{1*}$ to the sequence above,
we get the exact sequence
        $$0\to{\cal E}_{j-1}\to{\cal E}_j\to
        p_{1*}p_2^*\w_C((j+1)\Delta)|_{\Delta}\to 0.$$
Now, $p_\ell|_\Delta$ is an isomorphism for $\ell=1,2$. So
$p_{1*}p_2^*\w_C|_\Delta=\w_C$. In addition,
$p_{1*}{\cal O}_{C\times C}(-\Delta)|_\Delta=\w_C$. Thus
        \begin{eqnarray*}
    c_1({\cal E}_j)&=&c_1({\cal E}_{j-1})+
        c_1(p_{1*}p_2^*\w_C((j+1)\Delta)|_{\Delta})\\
    &=&-(j(j-1)/2)c_1(\w_C)+(1-(j+1))c_1(\w_C)\\
    &=&-(j(j+1)/2)c_1(\w_C),
    \end{eqnarray*}
as claimed.\qed

\begin{claim}\label{lem}{\bf Proposition.} 
Let $C$ be a general
smooth curve of genus $g\geq 1$ and $j$ a nonnegative
integer. Let $W_j\subseteq C\times C$ be the 
$j$-th Weierstrass divisor of $C$. Then $W_j$ 
is nonsingular and intersects 
the diagonal $\Delta$ transversally, at the pairs 
$(P,P)$ such that $P$ is a
Weierstrass point of $C$.
\end{claim}

\proof Let us show first that $W_j$ intersects 
$\Delta$ transversally. 
As pointed out in Subsection \ref{jWd}, the intersection 
$W_j\cap\Delta$ is,
set-theoretically, the set of pairs $(P,P)$ such that
$P$ is a Weierstrass point of $C$. As $C$ is general,
by Corollary~\ref{ws},
all its Weierstrass points are simple, and number
$g^3-g$ by Pl\"ucker Formula.
Now, since the intersection $W_j\cap\Delta$ is
finite, the number of points of intersection, weighted
by their intersection multiplicities, is equal to 
the degree
of the product $[W_j][\Delta]$. Using the notation 
and Formula (\ref{eq:17}) of Proposition \ref{prop7}, 
and using the Formulas
       \be
       \int_{C\times C}K_2[\Delta]=
       \int_{C\times C}K_1[\Delta]=
       -\int_{C\times C}[\Delta]^2=2g-2
       \label{eq:kd}
       \ee
and
       \be
       \int_{C\times C}K_1K_2=4(g-1)^2,
       \label{eq:kk}
       \ee
we get
       $$\int_{C\times C}[W_j][\Delta]=
       (g+j)(g+j+1)(g-1)-2j(g+j+1)(g-1)
        +j(j+1)(g-1),$$
which is exactly $g^3-g$. Thus the intersection
multiplicities are all one.

As a corollary of the transversal intersection, 
$W_j$ is nonsingular at its points on $\Delta$. So, 
let now $(P,Q)\in W_j$ for $P$ and $Q$ distinct, and 
let us show that $W_j$ is nonsingular at 
$(P,Q)$ as well. 

Let $J:={\rm Pic}^{g-1}(C)$, the component of the 
Picard scheme of $C$ parameterizing line bundles of 
degree $g-1$. Let $\Theta\subset J$ be the theta divisor, 
parameterizing line bundles with nontrivial global 
sections. Let
         $$\mu\: C\times C \to J$$
be the map taking a pair $(R,S)$ to the point of 
$J$ representing the bundle $\w_C((j+1)R-(g+j)S)$. 

We claim that $\mu(W_j)\subseteq\Theta$. Indeed, let 
$C^{(3)}:=C\times C\times C$, and denote by $p_{1,2}$ 
and $p_3$ the projection maps of 
$C^{(3)}$ onto the indicated factors. Let 
$\Delta_{1,3}$ and $\Delta_{2,3}$ be the indicated 
diagonals of $C^{(3)}$. Set
       $${\cal F}:=p_3^*\w_C((j+1)\Delta_{1,3}
         -(g+j)\Delta_{2,3}).$$
Recall the notation of Subsection \ref{wronski}. From 
the construction of $\Theta$, to show that 
$\mu(W_j)\subseteq\Theta$, it is enough to show that 
the Wronski map $\rho_{j,g+j-1}$ represents 
universally the cohomology of ${\cal F}$ or, put more 
simply, that $\rho_{j,g+j-1}$ can be viewed as a 
presentation of the right derived image 
$R^1(p_{1,2})_*{\cal F}$. 

Let ${\cal G}:=p_3^*\w_C((j+1)\Delta_{1,3})$. Then 
${\cal F}\subseteq{\cal G}$. From the definition of 
the Wronski map $\rho_{j,g+j-1}$, we get that 
$\rho_{j,g+j-1}$ is the image under $(p_{1,2})_*$ of 
the quotient map ${\cal G}\to{\cal G}/{\cal F}$. Thus, 
the map $\rho_{j,g+j-1}$ 
is the first map
in the following 
piece of the long derived sequence of 
$0\to{\cal F}\to{\cal G}\to{\cal G}/{\cal F}\to 0$ 
under $(p_{1,2})_*$:
        $$(p_{1,2})_*{\cal G}\to
          (p_{1,2})_*({\cal G}/{\cal F})\to
          R^1(p_{1,2})_*{\cal F}\to 
	  R^1(p_{1,2})_*{\cal G}.$$
Now, a fiberwise analysis shows that 
$R^1(p_{1,2})_*{\cal G}=0$. Thus $\rho_{j,g+j-1}$ 
is a presentation for $R^1(p_{1,2})_*{\cal F}$, finishing 
the proof that $\mu(W_j)\subseteq\Theta$. 

Let ${\cal L}:=\w_C((j+1)P-(g+j)Q)$, and denote 
by $[{\cal L}]$ the point of $J$ representing ${\cal L}$. 
Since $(P,Q)\in W_j$, we have $[{\cal L}]\in\Theta$. 
By Proposition \ref{WiPgen}, $h^0(C,{\cal L})=1$. Thus, 
it follows from \cite{acgh}, Prop. (4.2), p. 189, that 
$[{\cal L}]$ is a nonsingular point of $\Theta$. 
Furthermore, identifying the cotangent space of 
$J$ at $[{\cal L}]$ with $H^0(C,\w_C)$, the cotangent 
space of $\Theta$ at $[{\cal L}]$ is the quotient by 
the subspace $H^0(C,\w_C(-F))$, where $F$ is the 
unique effective divisor of $C$ such that 
${\cal L}={\cal O}_C(F)$. 

Identifying the cotangent space of $C\times C$ at 
$(P,Q)$ with $\w_C|_P\oplus\w_C|_Q$, the induced map 
of cotangent spaces 
$d\mu^*\: T^*_{J,[{\cal L}]}\to T^*_{C\times C,(P,Q)}$ 
is equivalent 
to the evaluation map,
       $$\epsilon\: H^0(C,\w_C)\to \w_C|_P\oplus\w_C|_Q.$$
We claim that $\epsilon(H^0(C,\w_C(-F)))\neq 0$. 
Indeed, if that were not the case, we would have 
$H^0(C,\w_C(-F))=H^0(C,\w_C(-F-P-Q))$, that is,
       $$h^0(C,{\cal O}_C((g+j-1)Q-(j+2)P))
         =h^0(C,{\cal O}_C((g+j)Q-(j+1)P)).$$
By the Riemann--Roch theorem,
       $$h^0(C,{\cal O}_C((g+j)Q-(j+1)P))=
         h^0(C,{\cal L})=1,$$
and thus, also by the Riemann--Roch theorem,
       $$h^0(C,\w_C((j+2)P-(g+j-1)Q))=3.$$
However, this contradicts Proposition \ref{arbD}. 

Since $\mu(W_j)\subseteq\Theta$, the image of 
$\epsilon(H^0(C,\w_C(-F)))$ in the cotangent space of 
$W_j$ at $(P,Q)$ is zero. Since 
$\epsilon(H^0(C,\w_C(-F)))\neq 0$, that cotangent space 
is a proper quotient of the cotangent space of 
$C\times C$ at $(P,Q)$, and thus has dimension at most 1. 
Since $W_j$ is a divisor, it follows that $W_j$ is 
nonsingular at $(P,Q)$.\qed

\section{Special ramification classes}\label{sec5}

\begin{claim}\label{srl} 
{\rm ({\it Special ramification loci.})
Let $C$ be a smooth curve of genus $g\geq 1$. 
For each nonnegative integer $i$, consider the 
following loci in $C\times C$:
\begin{enumerate}
\item The locus $D^+_i$ of pairs $(P,Q)\in C\times C$
such that
       $$(g+i-1)Q-(i+1)P$$
is linearly equivalent to an effective divisor.
\item The locus $E^+_i$ of pairs $(P,Q)\in C\times C$ such
that
       $$(g+i+1)Q-(i+1)P$$
is linearly equivalent to a moving effective divisor.
\item The locus $SW^+_i$ of pairs
$(P,Q)\in C\times C$ such 
that
$Q$ is a special ramification
point of $V_C(i,P)$.
\end{enumerate}
We claim that, set-theoretically,
       \be
       SW^+_i=D^+_i\cup E^+_i.\label{eq:swde}
       \ee

Indeed, by the Riemann--Roch Theorem, for a pair
$(P,Q)\in C\times C$, the divisor
$(g+i-1)Q-(i+1)P$ is linearly equivalent to an
effective one if and only if
        \be
    h^0(\w_C((i+1)P-(g+i-1)Q))\geq 2,\label{eq:6b3}
    \ee
while $(g+i+1)Q-(i+1)P$ is linearly equivalent to a
moving effective divisor if and only if
        \be
    h^0(\w_C((i+1)P-(g+i+1)Q))\geq 1.\label{eq:6b4}
    \ee
At any rate, if $(P,Q)\in D^+_i\cup E^+_i$, then $Q$ is a
special ramification point of $V_C(i,P)$, that is,
$(P,Q)\in SW^+_i$

On the other hand, let
$(P,Q)\in C\times C-(D^+_i\cup E^+_i)$. Then
        \begin{eqnarray*}
    &&h^0(\w_C((i+1)P-(g+i-1)Q))=1,\\
    &&h^0(\w_C((i+1)P-(g+i+1)Q))=0.
    \end{eqnarray*}
So, either $Q$ is an ordinary or a simple ramification
point of $V_C(i,P)$, that is, $(P,Q)\not\in SW^+_i$.

Let $\Delta$ be the diagonal subscheme of $C\times C$.
Notice that $E_i^+\cap\Delta$ consists of the pairs
$(P,P)$ such that $P$ is a Weierstrass point of $C$.
However, if $g>1$, both $D_i^+$ and $SW_i^+$ contain
$\Delta$. (If $g=1$, then 
$D_i^+=E_i^+=SW_i^+=\emptyset$.)

Let $D_i$, $E_i$ and $SW_i$ be the loci of points in
$D_i^+$, $E_i^+$ and $SW_i^+$ that lie off $\Delta$. Of
course, Expression (\ref{eq:swde}) implies
$SW_i=D_i\cup E_i$. Our Proposition~\ref{A2CC} claims
that, if $C$ is general, then
$SW_i=D_i\cup E_i$ holds in a more
refined way, in the cycle group of $C\times C$. Before
stating it, we need to endow $D_i$, $E_i$ and $SW_i$
with natural subscheme structures.}
\end{claim}

\begin{claim}\label{srs} 
{\rm ({\it Special ramification schemes}) 
Keep the notation of Subsection \ref{srl}, and recall 
that of Subsections \ref{wronski} and \ref{jWd}. 
Notice that the subsets 
$D_i^+$ and $E_i^+$ are the supports
of the degeneracy schemes of $\rho_{i,g+i-2}$ and
$\rho_{i,g+i}$, respectively. So we may give 
$D_i^+$ and $E_i^+$ the
corresponding scheme structures.
Give $D_i$ and $E_i$ the
corresponding open subscheme structures. We say that 
$D_i$ and $E_i$ are the {\it $i$-th special ramification 
schemes of type Diaz and Cukierman}, respectively. Call 
$E_i^+$ the {\it $i$-th expanded special ramification 
scheme of type Cukierman}.

In addition, differentiating along the fibers of $p_1$ a 
section of ${\cal O}_{C\times C}(Z_i)$ defining $Z_i$, 
we obtain a section of 
${\cal J}^1_{p_1}({\cal O}_{C\times C}(Z_i))$, 
well-defined modulo $\CC^*$. By functoriality, its zero 
scheme contains a pair $(P,Q)$ if and only if $Q$ is a 
special Weierstrass point of $V_C(i,P)$. Thus the zero 
scheme gives a scheme structure for
$SW_i^+$. Give $SW_i$ the induced open subscheme 
structure. We say that $SW_i$ is the {\it $i$-th special 
ramification scheme of $C$}.

Now, $Z_i=W_i+g\Delta$. As done 
for $Z_i$, we can differentiate along
the fibers of $p_1$ a section of ${\cal O}_{C\times C}(W_i)$
defining $W_i$ to obtain a section of ${\cal J}^1_{p_1}({\cal
O}_{C\times C}(W_i))$. Its zero scheme $S$ coincides with the
scheme $SW_i$ off $\Delta$, because $Z_i$ coincides with $W_i$
there. Moreover, if $C$ is general, then $S$ does not intersect
$\Delta$, and hence $S=SW_i$ scheme-theoretically. Indeed, let $P$
be a point of $C$. If $(P,P)\in W_i$, then $P$ is a Weierstrass
point of $C$. Moreover, as $C$ is general, by Corollary~\ref{ws},
the point $P$ is a simple Weierstrass point. So, it follows from
Proposition \ref{prop5} that $W_i$ intersects the fiber
$\{P\}\times C$ transversally at $(P,P)$. Thus the derivative
along $\{P\}\times C$ of a section defining $W_i$ does not vanish
at $(P,P)$. So $S\cap\Delta=\emptyset$.}
\end{claim}

\begin{claim}\label{matrices}{\bf Lemma.}
Let ${\cal O}$ be a local 
ring, and $r$ a nonnegative integer. Let $M$ be a matrix 
with $r+2$ rows and $r+1$ columns and entries in 
${\cal O}$. Let $M_1$ and $M_2$ be the submatrices 
obtained from $M$ by removing the last row, and
the last two rows, respectively. Assume that the 
matrix obtained from $M_1$ by taking residues 
has rank at least $r$. Let $z$
denote the determinant of $M_1$. Then there are 
$u,v\in{\cal O}$ such that
\begin{enumerate}
    \item $(z,u)$ is the ideal of all maximal minors 
      of $M_2$,
    \item $(z,v)$ is the ideal of all maximal minors 
      of $M$,
    \item $(z,uv)$ is the ideal generated by the two 
      maximal minors of $M$ obtained by removing each 
      of the last two rows.
\end{enumerate}
\end{claim}

\proof We may write $M$ in the form
        $$M=\left[\matrix{A&a&b\cr c&f_1&f_2
        \cr d&g_1&g_2\cr e&h_1&h_2}\right],$$
where $A$ is a square matrix of size $r-1$, where $a$ 
and $b$ are column vectors of 
size
$r-1$, where 
$c$, $d$ and $e$ are row vectors of 
size
$r-1$, 
and where $f_1$, $f_2$, $g_1$, $g_2$,
$h_1$ and $h_2$ are elements of ${\cal O}$.

Let $I$ and $J$ be the ideals of ${\cal O}$ generated,
respectively, by all maximal minors of the submatrices
        $$M_2=\left[\matrix{A&a&b\cr c&f_1&f_2}\right]
        \quad{\rm and }\quad M=\left[
        \matrix{A&a&b\cr c&f_1&f_2\cr d&g_1&g_2\cr
    e&h_1&h_2}\right].$$
Also, let $K\subseteq{\cal O}$ be the ideal generated by 
the determinants of the square submatrices
        $$M_1=\left[
        \matrix{A&a&b\cr c&f_1&f_2\cr d&g_1&g_2}
        \right]\quad{\rm and }\quad
        M'_1:=\left[\matrix
	{A&a&b\cr c&f_1&f_2\cr e&h_1&h_2}
        \right].$$
Notice that the determinant of the first matrix is $z$.

From the hypothesis, the matrix obtained from $M_2$ 
by taking residues has rank at least $r-1$. Thus, 
performing row and column operations on $M$, including 
column and row exchanges, we may, without changing the 
ideals $I$, $J$ and $K$, assume that $A$ is the 
identity matrix, $a=b=0$ and $c=d=e=0$. Then
$z=f_1g_2-f_2g_1$ and
    \begin{eqnarray*}
    I&=&(f_1,\, f_2),\\
    J&=&(f_1g_2-f_2g_1,\,f_1h_2-f_2h_1,\,g_1h_2-g_2h_1),\\
    K&=&(f_1g_2-f_2g_1,\, f_1h_2-f_2h_1).
    \end{eqnarray*}

Now, since the matrix obtained from $M_1$ 
by taking residues has rank at least $r$, at least one 
among $f_1,\, f_2,\, g_1,\, g_2$ is invertible.

If $f_1$ is invertible, then
        $$g_1h_2-g_2h_1=(g_1/f_1)(f_1h_2-f_2h_1)-
        (h_1/f_1)(f_1g_2-f_2g_1).$$
Thus, the lemma holds for $u=1$ and $v=f_1h_2-f_2h_1$. 
The case where $f_2$ is invertible is similar.

If $g_1$ is invertible, then
        \begin{eqnarray*}
        (f_1h_2-f_2h_1)&=&(f_1/g_1)(g_1h_2-g_2h_1)+
        (h_1/g_1)(f_1g_2-f_2g_1),\\
    f_2&=&(g_2/g_1)f_1-(1/g_1)(f_1g_2-f_2g_1).
        \end{eqnarray*}
Thus the lemma holds for $u=f_1$ and $v=g_1h_2-g_2h_1$. 
A similar analysis holds if $g_2$ is invertible.\qed

\begin{claim}\label{A2CC}{\bf Proposition.}
Let $C$ be a general
smooth curve of genus $g\geq 1$ and $i$ a nonnegative
integer. Let $\Delta$ be the diagonal of $C\times C$ and 
$W_i$ the $i$-th Weierstrass divisor. 
Let $SW_i$ be the $i$-th special ramification 
scheme, and $D_i$ and $E_i$ the $i$-th special 
ramification schemes of type Diaz and Cukierman, 
respectively. Let $E_i^+$ be the $i$-th expanded 
special ramification scheme of type Cukierman. 
Then these ramification schemes are finite and satisfy, 
in the cycle group of $C\times C$:
        $$[SW_i]=[D_i]+[E_i]\quad{\rm and }\quad
        [E_i^+]=[E_i]+(g+1)[W_i\cap\Delta].$$
\end{claim}

\proof Since $C$ is general, by Statement 1 of
Proposition~\ref{WiPgen}, the set $SW_i$ is finite
for each $i\geq 0$. Thus, so are $D_i$ and $E_i$
by Expression~(\ref{eq:swde}). It follows that
$E_i^+$ is finite, because $E_i^+\cap\Delta$ is the
set of points $(P,P)$ such that $P$ is Weierstrass, whence
is finite.

Recall the notation of Subsections \ref{wronski}, 
\ref{jWd}, \ref{srl} and \ref{srs}. 
Set $r:=g+i-1$. Both equalities can be proved locally. 
Thus, let $(P,Q)\in C\times C$ and ${\cal O}$ be the 
local ring of $C\times C$ at $(P,Q)$. 
As a map of ${\cal O}$-modules, 
$\rho_{i,r+1}$ is given by a matrix $M$ of the form 
described in the proof of Lemma \ref{matrices}. Let 
us use the notation described in the statement of 
that lemma.

Let $K\subseteq{\cal O}$ define $SW_i^+$. Then 
$K=(z,z')$, where $z$ (resp. $z'$) 
is the maximal minor obtained 
from $M$ by removing the last (resp. last but one) row. 
Notice that, 
from the nature of $M$ as a ``wronskian matrix'', $z'$ is 
also the derivative of $z$ along $p_1$. Let 
$I$ and $J$ be the ideals of ${\cal O}$ defining $D_i^+$ 
and $E_i^+$, respectively. Then $I$ and $J$ are the 
ideals of all the
maximal minors of $M_2$ and $M$, respectively.

Now, since $C$ is a general curve, by Statement 2 of 
Proposition \ref{WiPgen},
        $$h^0(\w_C((i+1)P-(g+i)Q))\leq 1.$$
This translates in the matrix obtained from $M_1$ by 
evaluating at $(P,Q)$ having rank at least $r$. 
Applying Lemma \ref{matrices},
there are $u,v\in{\cal O}$ such
    $$I=(z,u),\quad J=(z,v),\quad K=(z,uv).$$
Now, since $E_i^+$ is finite-dimensional and 
$C\times C$ is smooth, the sequence $z,v$ is regular. 
The same holds for the sequence $z,u$ 
if $P\neq Q$. It follows 
that $[SW_i]=[D_i]+[E_i]$.

The second equality in the statement of the
proposition is obvious off $\Delta$. Thus, assume 
$Q=P$. Since $E_i^+\cap\Delta=W_i\cap\Delta$, we may 
also assume that $P$ is a Weierstrass point of $C$. 

Let $t$ be a local parameter of $C$ at $P$, and 
$t_1,t_2\in{\cal O}$ be its pullbacks with respect to 
the projections $p_1$ and $p_2$. Then $t:=t_2-t_1$ is a 
local equation for $\Delta$. As we saw in Subsection 
\ref{srs}, we have $z=t^gw$, where
$w\in{\cal O}$ defines $W_i$, and is not divisible by 
$t$. Letting $\partial$ denote the 
derivative with respect to $t_2$, we have
       $$z'=\partial z=\partial(t^gw)=
       gt^{g-1}w+t^g\partial w.$$
Thus $t^{g-1}$ divides $z$ and $z'$, and hence each 
element of $K$, in particular $uv$. Since $E_i^+$ is 
finite, $t$ does not divide $v$, and hence $t^{g-1}|u$. 
Let $L:=t^{1-g}K$. Then there are two expressions for 
$L$:
       \be
       L=(tw,uv/t^{g-1})\quad{\rm and }\quad
       L=(tw,gw+t\partial w).
       \label{I}
       \ee
Since $W_i\cap\Delta$ is finite, the 
sequences $gw+t\partial w, t$ and $w, t$ are 
regular. Thus, from the second expression for $L$ above, 
we get
       $$\ell({\cal O}/L)=2\ell({\cal O}/(t,w))+
       \ell({\cal O}/(w,\partial w)).$$
Now, $\ell({\cal O}/(w,\partial w))=0$ because $w$ and 
$\partial w$ cut out $SW_i$, and $SW_i$ does not meet 
$\Delta$. Also, by Lemma \ref{lem}, $W_i$ intersects 
$\Delta$ transversally. Thus
$\ell({\cal O}/(t,w))=1$, and hence $\ell({\cal O}/L)=2$.

Now, since the sequence $z, v$ is regular, and 
$z=t^gw$, also the sequence $tw,v$ is regular. Thus, 
from the first expression for $L$ in (\ref{I}), 
we get
       $$\ell({\cal O}/L)=\ell({\cal O}/(tw,u/t^{g-1}))+
         \ell({\cal O}/(tw,v)),$$
and whence $\ell({\cal O}/(tw,v))\leq 2$. Since 
${\cal O}$ is regular, and the sequence $tw,v$ is 
regular, so is the sequence $v,w$. Thus
       $$\ell({\cal O}/(tw,v))=\ell({\cal O}/(t,v))+
         \ell({\cal O}/(w,v)).$$
Since $E_i^+$ contains $(P,P)$, the function $v$ is zero 
on $(P,P)$. Thus, since also $t$ and $w$ vanish on 
$(P,P)$, we get 
$\ell({\cal O}/(t,v))=\ell({\cal O}/(w,v))=1$. So, the 
multiplicity of $E_i^+$ at $(P,P)$ is
       $$\ell({\cal O}/(z,v))=
       g\ell({\cal O}/(t,v))+\ell({\cal O}/(w,v))
       =(g+1).$$
Since, by Lemma \ref{lem}, the multiplicity of 
$W_i\cap\Delta$ at $(P,P)$ is 1, we are done.\qed

\begin{claim}\label{prop07}{\bf Proposition.}
Let $C$ be a general
smooth curve of genus $g\geq 1$ and $i$ a nonnegative
integer. Let $SW_i$ be the $i$-th special ramification 
scheme of $C$. Then
        \be
    \int_{C\times C}[SW_i]=2ig(g-1)\Big((i+2)(g+i)^2
    +2(g+i)+2\Big).
    \label{eq:me1}
    \ee
\end{claim}

\proof Recall the notation of Subsections 
\ref{wronski}, \ref{jWd}, \ref{srl} and \ref{srs}. 
Since $C$ is general, $SW_i$ is finite. Also,
$SW_i$ is the zero scheme of a section of the
rank-2 bundle
${\cal J}^1_{p_1}({\cal O}_{C\times C}(W_i))$. Thus
its class in
the Chow group of $C\times C$ satisfies
        $$[SW_i]=c_2({\cal J}^1_{p_1}
        ({\cal O}_{C\times C}(W_i))).$$
Using the truncation sequence for bundles of jets, we get
        $$[SW_i]=[W_i](c_1(p_2^*\w_C)+[W_i]).$$
Now, $c_1(p_2^*\w_C)=K_2$. Using Expression
(\ref{eq:17}) for $j=i$, and taking into account
that $K_\ell^2=0$ for $\ell=1,2$, we get
        \begin{eqnarray*}
    [SW_i]&=&i(g+i+1)\Big((g+i)^2+g+i+1\Big)K_2[\Delta]\\
    &+&\frac{1}{2}i(i+1)\Big((g+i)^2+g+i+1\Big)K_1K_2\\
    &+&i^2(g+i+1)^2[\Delta]^2+i^2(g+i+1)(i+1)K_1[\Delta].
    \end{eqnarray*}
Using Formulas (\ref{eq:kd}) and (\ref{eq:kk}), we get
        \begin{eqnarray*}
    \int_{C\times C}[SW_i]&=&i(g+i+1)
    \Big((g+i)^2+g+i+1\Big)(2g-2)\\
    &+&\frac{1}{2}i(i+1)\Big((g+i)^2+g+i+1\Big)4(g-1)^2\\
    &-&i^2(g+i+1)^2(2g-2)+i^2(g+i+1)(i+1)(2g-2).
    \end{eqnarray*}
Simplifying, we get the claimed formula.\qed

\begin{claim}\label{thm6DiazCuk}{\bf Theorem.}
Let $C$ be a general smooth curve of genus $g\geq 1$,
and $i$ a nonnegative integer. Let $D_i$ and $E_i$ be the 
$i$-th special ramification schemes of type Diaz and 
Cukierman, respectively. Then $D_i$ and $E_i$ are 
reduced, and
      \be
      \int_{C\times C}[D_i]=g(g-1)
      \Big((g+i-1)^2(i+1)^2-(g-1)^2\Big)
      \label{eq:ex0}
      \ee
and
      \be
      \int_{C\times C}[E_i]=g(g-1)
      \Big((g+i+1)^2(i+1)^2-(g+1)^2\Big).
      \label{eq:ex6b0}
      \ee
\end{claim}

\proof Recall the notation of 
Subsections \ref{wronski}, \ref{jWd}, \ref{srl} and 
\ref{srs}. We will 
first compute the degrees of $D_i$ and $E_i$. 
First of all, since $E_i^+$ is finite, and
is the degeneracy scheme of $\rho_{i,g+i}$, applying
Porteous formula (\cite{Fu}, Thm. 14.4, p. 254), 
we get the
following expression for the class $[E_i^+]$ in the
Chow group of $C\times C$:
      $$[E_i^+]=c_2({\cal J}^{g+i}_{p_1}({\cal L}_i)-
      p_1^*{\cal E}_i).$$
Now, $c_2({\cal E}_i)=c_1({\cal E}_i)^2=0$,
since $C$ is one-dimensional. Thus
      $$[E_i^+]=c_2({\cal J}^{g+i}_{p_1}({\cal L}_i))
      -c_1({\cal J}^{g+i}_{p_1}({\cal L}_i))
      c_1(p_1^*{\cal E}_i).$$
Using the truncation sequence of the bundles of jets,
we get
      \begin{eqnarray*}
      c_1({\cal J}^{g+i}_{p_1}({\cal L}_i))&=&
      \sum_{\ell=1}^{g+i+1}(\ell K_2+(i+1)[\Delta]);\\
      c_2({\cal J}^{g+i}_{p_1}({\cal L}_i))&=&
      \sum_{m=2}^{g+i+1}\sum_{\ell=1}^{m-1}
      (\ell K_2+(i+1)[\Delta])(m K_2+(i+1)[\Delta]).
      \end{eqnarray*}
Expanding, and using that $K_2^2=0$, we get
      \begin{eqnarray*}
      c_1({\cal J}^{g+i}_{p_1}({\cal L}_i))&=&
      \frac{1}{2}(g+i+1)(g+i+2)K_2+(i+1)(g+i+1)[\Delta];\\
      c_2({\cal J}^{g+i}_{p_1}({\cal L}_i))&=&
      \frac{1}{2}(i+1)(g+i)(g+i+1)(g+i+2)K_2[\Delta]\\
      &&+\frac{1}{2}(i+1)^2(g+i)(g+i+1)[\Delta]^2.
      \end{eqnarray*}
Finally, using Formula (\ref{eq:c1Ei}) for $j=i$, and
Formulas (\ref{eq:kd}) and (\ref{eq:kk}), we get
      $$\int_{C\times C}[E_i^+]=(i+1)^2g(g-1)(g+i+1)^2.$$

Now, it follows from Proposition \ref{lem} that
$W_i$ meets $\Delta$
transversally at $g^3-g$ points. Thus, using
Proposition \ref{A2CC}, we get
      \begin{eqnarray*}
      \int_{C\times C}[E_i]&=&\int_{C\times C}[E_i^+]-
      (g+1)(g^3-g)\\
      &=&g(g-1)\Big((g+i+1)^2(i+1)^2-(g+1)^2\Big),
      \end{eqnarray*}
the stated formula for the degree of $[E_i]$.

Now, the expression for the degree of $[D_i]$ follows 
now from the equality $[SW_i]=[D_i]+[E_i]$ proved in
Proposition \ref{A2CC} and Formula (\ref{eq:me1})
for the degree of $[SW_i]$ proved in Proposition
\ref{prop07}.

Let us now show that $D_i$ and $E_i$ are reduced. 
Let $(P,Q)\in SW_i$. Let $\wh{\cal O}$ be the completion 
of the local ring of $C\times C$ at $(P,Q)$. 
Let $t_1$ and $t_2$ be local equations in $\wh{\cal O}$ 
for $\{P\}\times C$ and $C\times\{Q\}$, 
respectively. Then $\wh{\cal O}=\CC[[t_1,t_2]]$. 
Let $w\in{\cal O}$ be a local equation for 
$W_i$. Since $(P,Q)\in SW_i$, and since $W_i$ is 
nonsingular by Proposition \ref{lem}, we may assume that 
$w=t_1+u$, where $u\in\CC[[t_2]]$. Now, 
let $w'$ and $u'$ be the derivatives of $w$ and $u$ 
with respect to $t_2$. Then the ideal defining 
$SW_i$ at $(P,Q)$ is $(w,w')$, and the multiplicity of 
the cycle $[SW_i]$ at $(P,Q)$ is 
$\ell(\wh{\cal O}/(w,w'))$. Notice that $w'=u'$, and 
       $$\frac{\wh{\cal O}}{(w,w')}\cong
         \frac{\CC[[t_2]]}{(u')}=
         \frac{\CC[[t_2]]}{(u,u')}\cong
	 \frac{\CC[[t_1,t_2]]}{(t_1,w,w')}.$$
Thus the multiplicity of the cycle $[SW_i]$ at $(P,Q)$
is the multiplicity $m$ of $SW_i\cap(\{P\}\times C)$ at 
$(P,Q)$.

Since the formation of $SW_i$ commutes with base 
change, this multiplicity $m$ satisfies
       $$m=\wgt_V(Q)-1,$$
where $V$ is the complete linear system of sections 
of $\w_C((i+1)P)$. Now, by Propositions \ref{WiPgen} and 
\ref{arbD}, the order sequence of $V$ at $Q$ satisfies
       \begin{eqnarray*}
	 \epsilon_j(V,Q)&=&j\quad(j=0,1,\dots,g+i-3),\\
	 \epsilon_{g+i-2}(V,Q)&\leq&g+i-1,\\
	 \epsilon_{g+i-1}(V,Q)&\leq&g+i+1.
       \end{eqnarray*}
Thus $m\leq 2$, with equality if and only if 
       $$h^0(\w_C((i+1)P-(g+i-1)Q))=2\quad{\rm and}
         \quad h^0(\w_C((i+1)P-(g+i+1)Q))=1,$$
that is, if and only if $(P,Q)\in D_i\cap E_i$. 
Since $[SW_i]=[D_i]+[E_i]$ by Proposition~\ref{A2CC}, 
it follows that $D_i$ and $E_i$ are reduced at 
$(P,Q)$.\qed

\end{document}